
\magnification=1200
\overfullrule=0mm
\nopagenumbers
\vskip 1cm

\centerline {\bf Normal Systems of Algebraic and Partial
Differential Equations}
\bigskip
\centerline {\bf H. Hakopian}
\medskip
\centerline{Yerevan State University; King Saud University} 
\bigskip

In our talk we first discuss a joint paper with M. Tonoyan:
[On a Multivariate Theory, in "Approximation Theory: A volume dedicated to Blagovest Sendov", (B. Bojanov Ed.), Darba, Sofia, 2002, 212-230]. Here the polynomial interpolation approach is used to introduce the main results on multivariate normal algebraic systems. The
characterization of the algebraic systems which have maximal number of distinct and multiple solutions
is treated as the first and second versions of multivariate
fundamental theorem of algebra for normal systems, respectively. The relation between the second version and
a result of B. Mourrain is mentioned.

The main connection between normal algebraic systems 
and normal systems of PDEs is that the formers serve
as the systems of characteristic equations for the
latters (with constant coefficients). This multivariate
setting possesses all well-known properties of the classic
univariate case.
\medskip
Next we bring a construction which shows that any standard
algebraic system, with finite set of solutions,
$\ p(x,y)=0,\ q(x,y)=0,$ can be reduced to a normal
type algebraic system, namely to the following one
of total degree $n+m-1:$
$$
\psi_\alpha(x,y)p(x,y)=0,\phi_\beta(x,y)q(x,y)=0,\
\alpha,\beta\in Z_+^2, |\alpha|=n-1, |\beta|=m-1,
$$
where $n,m$ are the degrees of $p,q,$ respectively.
Here degrees of polynomials $\phi,\psi$ are $n-1,m-1$
and they have maximal number of zeros satisfying
the first and second equations of the initial system,
respectively. Beforehand we choose a coordinate system
such that the leading terms of $p$ and $q$ have no 
common (nontrivial) solution. This construction gives
some consequences interesting from the algebraic and
interpolation points of view.   
                                  
\bye